\documentclass[12pt,letterpaper]{amsart} %fix to 12pt

\usepackage{amsmath,amssymb,amsthm}
\usepackage[foot]{amsaddr}
\usepackage{fullpage}
\usepackage[mathscr]{eucal} % for \mathscr
\usepackage[noadjust]{cite} % to put citations in numerical order
\usepackage{hyperref}
\usepackage[capitalize, nameinlink]{cleveref}
\usepackage{enumitem}
\usepackage[usenames,dvipsnames]{xcolor} % for e.g. \todo color
\usepackage{tikz}
\usepackage[font=small,labelfont=bf, margin=.5in]{caption}

\colorlet{prettygreen}{ForestGreen!60!LimeGreen}

\usetikzlibrary{calc, arrows.meta}
\tikzset{vtx/.style={circle, fill, inner sep=1.5pt}}
\tikzset{openvtx/.style={circle, draw, inner sep=1.5pt}}

\usepackage{thmtools}

\newtheorem{theorem}{Theorem}[section]
\newtheorem{lemma}[theorem]{Lemma}
\newtheorem{proposition}[theorem]{Proposition}

\newtheorem{claim}[theorem]{Claim}
\newtheorem*{claim*}{Claim}

\theoremstyle{definition}
\newtheorem{definition}[theorem]{Definition}

\newtheorem{question}[theorem]{Question}

\newtheorem{conjecture}[theorem]{Conjecture}

\theoremstyle{remark}
\newtheorem{remark}[theorem]{Remark}

\crefname{claim}{Claim}{Claims}

\newcommand{\proofofclaim}{\renewcommand*{\qedsymbol}{$\blacksquare$}}

\newcommand{\eps}{\epsilon}
\newcommand{\LL}{{\mathcal L}}
\newcommand{\RP}{{\mathbb R\mathbb P}^2}
\newcommand{\TT}{{\mathbb T}^2}
\newcommand{\Vo}{V^\circ}

\newcommand{\skel}{{\mathcal S}}

\newcommand{\disk}{\mathbb D^2}

\DeclareMathOperator*{\E}{{\mathbb E}}
\DeclareMathOperator{\codeg}{codeg}
\DeclareMathOperator{\ex}{ex}

\newcommand{\exh}{\ex_{\hom}}

\title{An Improved Tur\'an Exponent for 2-Complexes}
\author{Maya Sankar}
\address{Department of Mathematics, Stanford University, Palo Alto, CA 94303.}
\email{mayars@stanford.edu}
\thanks{Sankar was supported by a Fannie and John Hertz Foundation Fellowship and NSF Graduate Research Fellowship DGE-1656518.}
\begin{document}
\maketitle

\begin{abstract}
The topological Tur\'an number $\ex_{\hom}(n,X)$ of a 2-dimensional simplicial complex $X$ asks for the maximum number of edges in an $n$-vertex 3-uniform hypergraph containing no triangulation of $X$ as a subgraph. We prove that the Tur\'an exponent of any such space $X$ is at most $8/3$, i.e., that $\ex_{\hom}(n,X)\leq Cn^{8/3}$ for some constant $C=C(X)$. This improves on the previous exponent of $3-1/5$, due to Keevash, Long, Narayanan, and Scott. Additionally, we present new streamlined proofs of the asymptotically tight upper bounds for the topological Tur\'an numbers of the torus and real projective plane, which can be used to derive asymptotically tight upper bounds for all surfaces. 
The key insight is an improved understanding of the placement of 4-cycles $vwv'w'$ that are likely to bound a triangulation of the disk within a randomly-selected subset of vertices.
\end{abstract}

\section{Introduction}

Given a graph $F$, the standard Tur\'an problem asks for the maximum number of edges in an $n$-vertex graph containing no subgraph isomorphic to $F$. Originally studied by Tur\'an \cite{Tu41} in 1941, this problem has become a cornerstone of extremal combinatorics. A seminal result of Erd\H os, Stone, and Simonovits \cite{ErSi66,ErSt46} determines this quantity asymptotically in terms of the chromatic number $\chi(F)$ for all graphs with chromatic number at least 3. For bipartite $F$, the so-called degenerate Tur\'an problem remains an active area of study (see e.g.\ \cite{JaSu24,JJM22,SuTo20}).

Tur\'an's question generalizes naturally to hypergraphs. Given a fixed $r$-uniform hypergraph $F$, its \emph{Tur\'an number} $\ex(n,F)$ is the maximum number of edges in an $r$-uniform hypergraph on $n$ vertices with no subgraph isomorphic to $F$. Tur\'an numbers of hypergraphs form a rich and fruitful area of study; for more background on this topic, we refer the reader to the survey of Keevash \cite{Ke11}. This paper is concerned with a variant of the standard Tur\'an problem concerning topological spaces, whose analysis requires a blend of topological and probabilistic ideas.

One may view an $r$-uniform hypergraph $H$ as an $(r-1)$-dimensional simplicial complex on the same vertex set whose facets are the edges of $H$ --- i.e., any set of vertices contained in an $r$-edge of $H$ forms a simplex.
Given a fixed $(r-1)$-dimensional simplicial complex $X$, a \emph{homeomorph} of $X$ is a subgraph of $H$ that is homeomorphic (as a simplicial complex) to $X$. It is natural to consider the following extremal quantity. Let $\exh(n,X)$ denote the maximum number of edges in an $r$-uniform hypergraph $H$ containing no homeomorph of $X$; this is the \emph{topological Tur\'an number} of the topological space $X$.
This quantity can also be viewed as an ordinary Tur\'an number: if we write $\mathcal F_X$ for the family of $r$-uniform hypergraphs homeomorphic to $X$, then we have $\exh(n,X)=\ex(n,\mathcal F_X)$. Alternatively, from the topologist's perspective, it is most natural to interpret $\exh(n,X)$ as the maximum number of $(r-1)$-dimensional faces in an $n$-vertex simplicial complex with no subcomplex homeomorphic to $X$.

In the 1-dimensional case, $\exh(n,X)$ is determined asymptotically by an old result of Mader \cite{Ma67}. He proved that if a graph $G$ contains no subdivision of $K_t$ then $G$ has at most $cn$ edges, where the constant $c=c(t)$ depends only on $t$. In the language of topological Tur\'an numbers, we write $\exh(n,K_t)=\Theta(n)$, where $K_t$ is viewed as a 1-dimensional simplicial complex (a space prwhich is, up to homotopy, a wedge of $\binom {t-1}2$ copies of $\mathbb S^1$).

The study of $\exh(n,X)$ in higher dimensions was initiated in 1973 by S\'os, Erd\H os, and Brown \cite{SoErBr73}, who showed that the topological Tur\'an number of the sphere $\mathbb S^2$ is asymptotically $\Theta(n^{5/2})$. Motivated by this result, Linial \cite{Li08,Li18} repeatedly asked for upper bounds on $\exh(n,X)$ for other spaces $X$, in particular the torus $\TT$, as part of his program in high-dimensional combinatorics. This was recently resolved asymptotically for surfaces by work of Kupavskii, Polyanskii, Tomon, and Zakharov \cite{KPTZ21} and the author \cite{Sa22}, who showed that $\exh(n,X)=\Theta(n^{5/2})$ for all orientable and non-orientable surfaces, respectively.

This paper focuses on the topological Tur\'an exponent of an arbitrary 2-dimensional simplicial complex $X$. We remark that Mader's result proves a ``universal'' exponent in the 1-dimensional case: the dependency of $\exh(n,X)$ on $n$ is linear for any 1-dimensional simplicial complex $X$.
In 2 dimensions, Keevash, Long, Narayanan, and Scott \cite{KLNS21} proved that $\exh(n,X)=O(n^{3-1/5})$ for any fixed 2-dimensional simplicial complex $X$. 

Our main result improves this exponent to 8/3. Let $K_t^{(3)}$ denote the complete 3-uniform hypergraph on $t$ vertices (or its associated 2-dimensional simplicial complex).
\begin{theorem}\label{mainthm}
For any $t$, there is a constant $C=C(t)$ such that $\exh\left(n,K_t^{(3)}\right)\leq Cn^{8/3}$.
\end{theorem}

Let us briefly mention what is known in higher dimensions. Long, Narayanan, and Yap \cite{LNY22} recently proved the existence of universal exponents in every dimension. That is, they showed that every $d$-dimensional simplicial complex has topological Tur\'an exponent at most $d+1-\eps_d$, for some constant $\eps_d$ depending only on the dimension $d$. They attained the lower bound $\eps_d\geq d^{-2d^2}$; however, it is unclear what value of $\eps_d$ might be optimal. Even the asymptotics of $\exh(n,\mathbb S^d)$ remain unknown, though recent work of Newman and Pavelka \cite{NePa24} suggests that the correct exponent lies between $d+1-(d+1)/(2^{d+1}-2)$ and $d+1-1/2^{d-1}$.

We remark that the best possible exponent in \cref{mainthm} is 5/2, matching the bound for spheres and other surfaces. Reiterating a conjecture of Keevash, Long, Narayanan, and Scott, we suspect that this is the correct exponent. 
\begin{conjecture}[{\cite[Conjecture 1.2]{KLNS21}}]\label{conj:upperbound}
	Let $X$ be any 2-dimensional simplicial complex. There is a constant $C=C(X)$ such that $\exh(n,X)\leq Cn^{5/2}$.
\end{conjecture}

\noindent
This conjecture is supported in the random setting: Gundert and Wagner \cite{GW16} showed that the threshold for the random 3-uniform hypergraph $G^{(3)}(n,p)$ to contain a homeomorph of $K_t^{(3)}$ is $p=\Theta(1/\sqrt n)$, where the constant coefficient depends on $t$.

Our proof of \cref{mainthm} introduces some new ideas that we feel may be useful to further improve the universal exponent in the 2-dimensional setting. Both \cite{KLNS21} and \cite{LNY22} obtain universal exponents by defining a specific $r$-uniform $r$-partite hypergraph $F_t^{(r)}$ homeomorphic to $K_t^{(r)}$ and locating a subgraph isomorphic to $F_t^{(r)}$. 
In contrast, our construction glues together homeomorphs of the disk $\disk$ whose isomorphism classes (as hypergraphs) remain unspecified. The heart of our method is an improved understanding of \emph{disk-coverable cycles}, which were tacitly studied in \cite{KPTZ21,Sa22}. Defined formally in \cref{sec:prelim}, such a cycle is the boundary of many different homeomorphs of the disk, which may not be isomorphic as hypergraphs. Our approach is to locate a suitable configuration of disk-coverable cycles, whereupon we build a homeomorph of $K_t^{(3)}$ as a union of disk homeomorphs bounded by the disk-coverable cycles.

As another application of our analysis of disk-coverable cycles, we derive shorter proofs of the upper bounds for the topological Tur\'an numbers of the torus $\TT$ and real projective plane $\RP$. These two bounds are the main results of \cite{KPTZ21} and \cite{Sa22}, and can be leveraged to prove an exponent of $5/2$ for all orientable and nonorientable surfaces, respectively.
\begin{theorem}[\cite{KPTZ21} and \cite{Sa22}]\label{thm:surfaces}
	We have $\exh(n,\TT)=O(n^{5/2})$ and $\exh(n,\RP)=O(n^{5/2}$).
\end{theorem}
\noindent

The remainder of this paper is organized as follows. In \cref{sec:prelim}, we give some preliminary definitions, allowing us to outline our approach more formally. \cref{sec:diskc} develops our understanding of disk-coverability, proving our key technical result. Using these preliminaries, we prove \cref{mainthm} in \cref{sec:mainpf}. \cref{sec:TandRP} contains our new proof of \cref{thm:surfaces}, which also leverages the key result of \cref{sec:diskc}. Lastly, \cref{sec:end} concludes with a discussion of the limitations of our approach --- highlighting connections to the study of bipartite Tur\'an exponents --- and some open problems.

\section{Preliminaries}\label{sec:prelim}

In this section, we sketch the main ideas behind our result and introduce some key definitions.

Let $H$ be a 3-uniform hypergraph. Its \emph{1-skeleton} $\skel(H)$ is the graph on $V(H)$ with edge set
\[
E(\skel(H))=\{xy:xyz\in E(H)\text{ for some }z\in V(H)\}.
\]
For a vertex $u\in V(H)$, its \emph{link graph} $H_u$ is the subgraph of $\skel(H)$ on vertex set $V(H)\setminus\{u\}$ whose edges $vw\in E(H_u)$ correspond to edges $uvw\in E(H)$ containing $u$. 

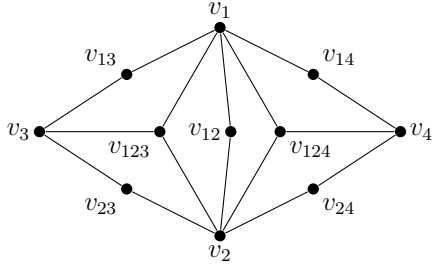
\begin{figure}
\begin{tikzpicture}[scale=1.6]
	\path (0,0) coordinate[vtx] (v123)% node[right] {$v_{123}$}
		--+ (60:1) coordinate[vtx] (v1)
		--+ (-60:1) coordinate[vtx] (v2)
		--+ (180:1) coordinate[vtx] (v3)
		--+ (120:0.55) coordinate[vtx] (v13)% node[above left] {$v_{13}$}
		--+ (-120:0.55) coordinate[vtx] (v23)% node[below left] {$v_{23}$}
		--+ (0:0.59) coordinate[vtx] (v12)
		--++(0:1) coordinate[vtx] (v124)
		--+ (0:1) coordinate[vtx] (v4)
		--+ (60:0.55) coordinate[vtx] (v14)
		--+ (-60:0.55) coordinate[vtx] (v24);
	\draw (v1) -- (v12) -- (v2) -- (v23) -- (v3) -- (v13) -- (v1)
		-- (v14) -- (v4) -- (v24) -- (v2);
	\draw (v123) -- (v1) (v123) -- (v2) (v123) -- (v3);
	\draw (v124) -- (v1) (v124) -- (v2) (v124) -- (v4);
	\foreach \idx/\pos in {
		1/above,2/below,3/left,4/right,
		12/left,
		13/above left, 23/below left, 14/above right, 24/below right,
		123/below left, 124/below right} {
		\node[\pos, scale=0.85] at (v\idx) {$v_{\idx}$};
	}
\end{tikzpicture}
\caption{Part of $\Gamma_4$}
\label{fig:Gamma}
\end{figure}

Our proof relies on the following topological decomposition of $K_t^{(3)}$.
Build a graph $\Gamma_t$ with $\binom t1 + \binom t2 + \binom t3$ vertices. Its vertices are of the form $v_i$, $v_{ij}$, or $v_{ijk}$ with $i,j,k\in[t]$ distinct, corresponding to the 0-, 1-, and 2-simplices of $K_t^{(3)}$ respectively. Connect each vertex $v_{ij}$ to $v_i$ and $v_j$ and connect each $v_{ijk}$ to $v_i$, $v_j$, and $v_k$, as pictured in \cref{fig:Gamma}. We note that $\Gamma_t$ is bipartite, with the $t$ vertices $v_1,\ldots,v_t$ forming one side of the bipartition. 
One may build $K_t^{(3)}$ as a CW complex with 1-skeleton $\Gamma_t$ by attaching one disk to each of the $3\binom r3$ cycles of the form $v_iv_{ij}v_jv_{ijk}$. Call these $3\binom r3$ cycles the \emph{special 4-cycles} of $\Gamma_t$.

To locate a homeomorphic copy of $K_t^{(3)}$ in $H$, we find a copy of $\Gamma_t$ in $\skel(H)$ --- more specifically, we embed $\Gamma_t$ in a carefully chosen link $H_u$ --- such that each special 4-cycle $C\subseteq\Gamma_t$ bounds many disks $D\subseteq H$. We formalize this notion (called disk-coverability) in the next few paragraphs.

In the hypergraph setting, a \emph{disk} $D$ is a 3-uniform hypergraph homeomorphic to the (topological) disk $\disk$. Its \emph{boundary}, which is a cycle in $\skel(D)$, is denoted $\partial D$, and its \emph{interior vertex set} is $\Vo(D)=V(D)\setminus V(\partial D)$. For technical reasons, we introduce one further property. A disk $D$ with at least two edges is said to be \emph{boundary-inducing} if $\partial D$ is an induced subgraph of $\skel(D)$; equivalently, $\partial D$ is an induced subcomplex of the simplicial complex associated to $D$.

We have now built up the requisite terminology to define disk-coverability. Informally, this is the property that a cycle in $\skel(H)$ bounds many disks in $H$.

\begin{definition}\label{def:diskc}
Fix $p,\eps\in(0,1)$ and let $H$ be a 3-uniform hypergraph. Choose $U\subseteq V(H)$ by selecting each vertex independently with probability $p$. For a cycle $C\subseteq\skel(H)$, consider the event $A_C$ that there is a boundary-inducing disk $D\subseteq H$ with $\partial D=C$ and $\Vo(D)\subseteq U$. We say the cycle $C$ is \emph{$(p,\eps)$-disk-coverable} if $\Pr[A_C]\geq 1-\eps$.
\end{definition}

Formally, we prove \cref{mainthm} by locating a copy of $\Gamma_t$ in $\skel(H)$ such that each special 4-cycle of $\Gamma_t$ is $(p,\eps)$-disk-coverable, where the parameters $p=t^{-3}$ and $\eps=2t^{-6}$ depend on $t$. The proof concludes by randomly partitioning $V(H)$ into $1/p=t^3$ parts $U_1,\ldots,U_{t^3}$ and identifying $3\binom t3$ disks bounded by the special 4-cycles of $\Gamma_t$, each with interior vertices in a different part $U_i$.

We briefly comment on the differences between our work and Keevash, Long, Narayanan, and Scott's work deriving a universal exponent of $3-1/5$ \cite{KLNS21}. The initial setup is nearly identical: they too find a copy of $\Gamma_t$ in $\skel(H)$ and then glue a disk to each special 4-cycle. However, they only consider disks formed by four triangles sharing a common vertex. Hence, their embedding of $\Gamma_t$ requires that every special 4-cycle is contained in many different link graphs $H_u$, inherently requiring that $H$ contains many copies of the octahedron $K_{2,2,2}^{(3)}$. We surpass this limitation by considering a broader class of disks. Indeed, it is conjectured that the extremal number of the octahedron is $\Theta(n^{3-1/4})$, although the best known lower bound is $\Omega(n^{8/3})$ due to Katz, Krop, and Maggioni \cite{KKM02}, which matches our upper bound.

\section{Disk-Coverability and Admissibility}\label{sec:diskc}

In this section, we collect technical results regarding disk-coverability that prove useful in \cref{sec:mainpf,sec:TandRP}. We first present \cref{lem:embedding}, which motivates the study of disk-coverability. This lemma allows us to build homeomorphs by attaching disks to a 1-skeleton. Then, we relate disk-coverability to \emph{admissibility}, a notion introduced by Kupavskii, Polyanskii, Tomon, and Zakharov in \cite{KPTZ21}. Our key result concerning admissibility is \cref{lem:2adm}, which shows that inadmissible structures are relatively rare. In \cref{sec:mainpf}, this allows us to bound the number of non-disk-coverable 4-cycles in the average link of a 3-uniform hypergraph.

One should think of \cref{lem:embedding} as formalizing the last step of our proof sketch. At this point, we have embedded $\Gamma_t$ in the 1-skeleton of a 3-uniform hypergraph such that each special 4-cycle is $(p,\eps)$-disk-coverable. \cref{lem:embedding} shows that one may attach disks to all of the special 4-cycles in such a way that their union forms a homeomorph of $K_t^{(3)}$. While this seems intuitive, one must verify that there are no spurious intersections between the disks' interiors.

\begin{lemma}\label{lem:embedding}
Let $H$ be a 3-uniform hypergraph, and fix cycles $C_1,\ldots,C_k\subseteq\skel(H)$ in its 1-skeleton. Set $m=|\bigcup_{i=1}^kV(C_i)|$, and suppose there are fixed parameters $p\leq\frac 1{k+m}$ and $\eps<\frac 1{k(k+m)}$ such that each cycle $C_i$ is $(p,\eps)$-disk-coverable. Then there are disks $D_1,\ldots,D_k$ with $\partial D_i=C_i$, such that for any $i\neq j$, the simplicial complexes associated to $D_i$ and $D_j$ only intersect on the one-dimensional simplicial complex $C_i\cap C_j$.
\end{lemma}

\begin{proof}
Partition $V(H)$ into $k+m$ sets $U_1,\ldots,U_{k+m}$ by independently placing each $v\in V(H)$ into a part $U_i$ chosen uniformly at random. For $i\in[k]$ and $j\in[k+m]$, let $A_{i,j}$ be the event that there is a boundary-inducing disk $D_{i,j}$ with $\partial D_{i,j}=C_i$ and $\Vo(D)\subseteq U_j$. Because $C_i$ is $(p,\eps)$-disk-coverable with $p\leq\frac 1{k+m}$, it follows that $\Pr[A_{i,j}]\geq 1-\eps$. A union bound yields
\[
\Pr[A_{i,j}\text{ hold simultaneously for all }i,j]\geq 1-k(k+m)\eps>0.
\]
Thus, there is a partition $V(H)=U_1\cup\cdots\cup U_{k+m}$ such that the disks $D_{i,j}$ exist for each $i\in[k]$ and $j\in[k+m]$.

Let $W=\bigcup_{i=1}^kV(C_i)$. Because $W$ has size $m$, we may reorder the $k+m$ sets $U_i$ such that $W$ does not intersect $U_1,\ldots,U_k$. Set $D_i=D_{i,i}$, so that $D_i$ is a boundary-inducing disk with $\partial D_i=C_i$ and $\Vo(D_i)\subseteq U_i$.

We verify that the simplicial complexes corresponding to these disks cannot intersect except on their boundaries. Fix $i\neq j$ and let $X_i$ and $X_j$ be the simplicial complexes associated to the disks $D_i$ and $D_j$.
We note that $\Vo(D_i)\subseteq U_i$ is disjoint from $V(X_j)=\Vo(D_j)\cup V(\partial D_j)\subseteq U_j\cup W$. Moreover, the restriction of $X_i$ to the vertex set $V(X_i)\setminus\Vo(D_i)=V(\partial D_i)$ is exactly $\partial D_i=C_i$, because $D_i$ is a boundary-inducing disk. It follows that $X_i\cap X_j$ is a subcomplex of the 1-dimensional simplicial complex $C_i$. Similarly, $X_i\cap X_j$ is a subcomplex of $C_j$, and it follows that the intersection $X_i\cap X_j$ is exactly $C_i\cap C_j$.
\end{proof}

The remainder of this section develops techniques to identify disk-coverable 4-cycles, building off of the following insight of Kupavskii, Polyanskii, Tomon, and Zakharov \cite{KPTZ21}. Suppose $v$ and $v'$ are vertices of a 3-uniform hypergraph $H$. Let $G=H_v\cap H_{v'}$ and suppose vertices $w,w'\in V(G)$ are connected by many paths in $G$. Then there are many disks with boundary $vwv'w'$: for any path $w_0\cdots w_k$ from $w_0=w$ to $w_k=w'$, the disk $D$ induced by the edge set
\begin{equation}\label{eqn:pyramid}
E(D)=\bigcup_{i=0}^{k-1}\{vw_iw_{i+1}, v'w_iw_{i+1}\}	
\end{equation}
has boundary $\partial D=vw_0v'w_k=vwv'w'$. In \cite{KPTZ21}, the authors showed that most edges $ww'\in E(G)$ have many paths from $w$ to $w'$, implying that that the 4-cycle $vwv'w'\subset\skel(H)$ is disk-coverable for most edges $ww'\in E(G)$. They called such an edge \emph{admissible} in $G$.

We extend this idea to paths of length 2 in $G=H_v\cap H_{v'}$. Informally, a $P_2$ in $G$ is \emph{admissible} (see \cref{def:adm} below) if it is contained in many cycles of length at least 4. Controlling the number of inadmissible $P_2$'s in $G$ is substantially more challenging than bounding the number of inadmissible edges --- indeed, as discussed later in this section, there are graphs $G$ in which almost all $P_2$'s are inadmissible. In \cref{lem:2adm}, we bypass this impediment by counting $P_2$'s according to a suitable weighting factor. Finding admissible paths of length 2 is vital to our proof of \cref{mainthm} in \cref{sec:mainpf} and is the key ingredient that allows us to streamline the proof of \cref{thm:surfaces} in \cref{sec:TandRP}.

The remainder of this section is dedicated to stating and proving \cref{lem:2adm}.

\begin{definition}\label{def:adm}
	Fix $p,\eps\in(0,1)$. Let $G$ be a graph and $wuw'$ a path of length 2 in $G$. Choose $U\subseteq V(G)\setminus\{u\}$ by selecting each vertex except $u$ independently with probability $p$. Consider the event $A_{wuw'}$ that $G$ contains a path $ww_1\cdots w_{k-1}w'$ of length $k\geq 2$ from $w$ to $w'$ such that $w_1,\ldots,w_{k-1}\in U$. We say the path $wuw'$ is \emph{$(p,\eps)$-admissible} if $\Pr[A_{wuw'}]\geq 1-\eps$.
\end{definition}

As observed earlier, any path from $w$ to $w'$ in $G=H_v\cap H_{v'}$ gives rise to a disk with edge set given by (\ref{eqn:pyramid}) and boundary $vwv'w'$. It is straightforward to check that this disk is boundary-inducing if the path from $w$ to $w'$ has length at least 2. This yields the following relationship between admissibility and disk-coverability.

\begin{proposition}\label{prop:diskc}
	Let $H$ be a 3-graph and fix $v,v'\in V(H)$. Suppose a path $wuw'$ is $(p,\eps)$-admissible in $H_v\cap H_{v'}$. Then the 4-cycle $vwv'w'$ is $(p,\eps)$-disk-coverable in $H$.
\end{proposition}

Having defined admissibility, we may now state \cref{lem:2adm}. In the following, the notation
\[
\sum_{xyz\in A}f(xyz)
\]
refers to a sum taken over \emph{unlabeled} paths of length 2 in some set $A$. That is, each term $f(xyz)=f(zyx)$ occurs exactly once in the sum.

\begin{lemma}\label{lem:2adm}
Fix $p,\eps\in(0,1)$ and let $G$ be a graph on $n$ vertices. Let $A$ be the set of paths of length 2 in $G$ that are not $(p,\eps)$-admissible. Then
\[
\sum_{xyz\in A}\frac 1{\deg y}<\frac {3n}{2p^2\eps}.
\]
\end{lemma}

Before proving the lemma, let us briefly discuss the $\frac 1{\deg y}$ term. 
We note that it is impossible to derive a useful bound without this weighting factor. In fact, there are $n$-vertex graphs with a linear number of edges in which a $(1-o(1))$-fraction of $P_2$'s are inadmissible for any parameters $p,\eps$. For example, start with a clique on $\sqrt n$ vertices and attach $n-\sqrt n$ pendant vertices to one vertex of the clique. A $(1-o(1))$-fraction of the $P_2$'s in this graph are contained in no cycle, and are thus inadmissible for any parameters $p,\eps$.

We first prove \cref{lem:2adm} in the limiting case that $p=1$ and $\eps$ approaches 1. The choice of weight $\frac 1{\deg y}$ in \cref{lem:2adm} is motivated by this proof.

\begin{proposition}[\cref{lem:2adm} with $p=1$, $\eps\to 1^-$]\label{prop:P2s}
Let $G$ be a graph with $n$ vertices. Let $A$ be the set of paths of length 2 in $G$ that are not contained in a cycle of length at least 4. Then,
\[
\sum_{xyz\in A}\frac 1{\deg y}<\frac{3n}2.
\]
\end{proposition}

\begin{proof}
We may assume that $G$ is connected --- if not, the following proof can be applied separately to each connected component of $G$.

Fix a vertex $v_0$ of $G$. Let $T$ be a spanning tree of $G$ rooted at $v_0$ with the following property: the vertices at depth $d$ in $T$ are exactly those vertices at distance $d$ from $v_0$ in $G$. Such a $T$ can be constructed by, for each vertex $v$ at distance $d>0$ from $v_0$, choosing its parent among the neighbors of $v$ at distance $d-1$ from $v_0$. Write $\gamma(v)$ for the set of children of $v$ in $T$ and $\gamma^*(v)\supseteq\gamma(v)$ for the set of descendants of $v$ (not including $v$ itself). In the following, $N(v)$ and $\deg(v)$ always denote neighborhoods and degrees computed in the original graph $G$, not in the spanning tree $T$.

Consider a path $xyz\in A$. Because there is no 4-cycle in $G$ containing $xyz$, it follows that $y$ is the unique common neighbor of $x$ and $z$. We additionally claim that $xz\in E(T)$ if neither $x$ nor $z$ is in $\gamma^*(y)$. Observe that $T\setminus\{y\}$ has a connected component with vertex set $V(G)\setminus(\{y\}\cup \gamma^*(y))$. Hence, if $x,z\notin\gamma^*(y)$, then $x$ and $z$ are connected by a path in $T\setminus\{y\}$ of length $\ell\geq 1$. By adding the edges $xy$ and $yz$ to this path, one obtains a cycle of length $\ell+2$ containing $xyz$; however, there is no such cycle of length at least 4. We conclude that $\ell=1$, i.e., that $xz$ is an edge of $T$.

From the prior paragraph, it follows that each path $xyz\in A$ either has an endpoint in $\gamma^*(y)$ or satisfies $xz\in E(T)$ with $y$ the unique common neighbor of $x$ and $z$. Thus,
\begin{align*}
\sum_{xyz\in A}\frac 1{\deg y}
&\leq\sum_{y\in V(G)}\frac{
	|\gamma^*(y)\cap N(y)|\times(\deg y) + |\{xz\in E(T):N(x)\cap N(z)=\{y\}\}|
}
{\deg y}
\\&\leq\left(\sum_{y\in V(G)}|\gamma^*(y)\cap N(y)|\right) + \frac {|\{xz\in E(T):|N(x)\cap N(z)|=1\}|}2,
\end{align*}
where we bound the second term using $\deg y\geq 2$.
To evaluate this sum, we show that $\gamma^*(y)\cap N(y)\subseteq\gamma(y)$. (These sets are in fact equal, but we do not need the opposite containment.) Suppose that $y$ is at depth $d$ in $T$, i.e., that $y$ is at distance $d$ from $v_0$ in $G$. Then every vertex of $N(y)$ has distance at most $d+1$ from $v_0$ in $G$, and hence is at depth at most $d+1$ in $T$. Conversely, every vertex of $\gamma^*(y)\setminus \gamma(y)$ has depth at least $d+2$ in $T$. It follows that the sets $N(y)$ and $\gamma^*(y)\setminus\gamma(y)$ are disjoint, implying that $\gamma^*(y)\cap N(y)\subseteq\gamma(y)$.

Therefore,
\[
\sum_{xyz\in A}\frac 1{\deg y}
\leq\left(\sum_{y\in V(G)}|\gamma(y)|\right)+ \frac{|\{xz\in E(T):|N(x)\cap N(z)|=1\}|}2
<n+\frac{|E(T)|}2<\frac{3n}2.
\qedhere
\]
\end{proof}

We now leverage \cref{prop:P2s} to prove \cref{lem:2adm} for any choice of parameters $p,\eps$.

\begin{proof}[Proof of \cref{lem:2adm}]
Choose $U\subseteq V(G)$ by selecting each vertex independently with probability $p$. Let $A_U$ be the set of paths of length 2 in the induced subgraph $G[U]$ that are not contained in a cycle of length at least 4 in $G[U]$. \cref{prop:P2s} implies that
\[
\frac{3|U|}2>\sum_{xyz\in A_U}\frac 1{\deg_Uy}\geq\sum_{xyz\in A_U}\frac 1{\deg_Gy}.
\]
%Recall that $A$ is the set of paths of length 2 that are not $(p,\eps)$-admissible in $G$. Thus, 
Taking the expectation of both sides yields
\begin{align*}
\frac{3pn}2
&>\sum_{\substack{xyz\subset G\\\text{a }P_2}}\frac{\Pr[xyz\in A_U]}{\deg_Gy}
\geq\sum_{xyz\in A}\frac{\Pr[xyz\in A_U]}{\deg_Gy}
\\&=\sum_{xyz\in A}
	\frac{\Pr[x,y,z\in U\text{ and $\nexists$ cycle $C\subset G[U]$ with $xyz\subset C$ and $|C|\geq 4$}]}
	{\deg_G y}
\\&=\sum_{xyz\in A}
	\frac{\Pr[x,y,z\in U]\times\Pr[\text{$\nexists$ $C\subset G[U\cup\{x,y,z\}]$ with $xyz\subset C$ and $|C|\geq 4$}]}
	{\deg_Gy}.
%\\&\geq\sum_{xyz\in A}\frac{p^3\times\eps}{\deg_Gy}.
\end{align*}
To bound the numerator, recall that $A$ is the set of non-$(p,\eps)$-admissible $xyz\subset G$. In other words, for each $xyz\in A$, there exists such a cycle $C\subseteq G[U\cup\{x,y,z\}]$ with probability at most $1-\eps$. Thus, the numerator is at least $p^3\times \eps$ and
\[
\frac{3pn}{2}>\sum_{xyz\in A}\frac{p^3\times\eps}{\deg_Gy},
\]
which rearranges to
\[
\sum_{xyz\in A}\frac{1}{\deg_Gy}<\frac{3n}{2p^2\eps}.\qedhere
\]
\end{proof}

\section{Proof of the Main Result}\label{sec:mainpf}

In this section, we prove \cref{mainthm}. Let $H$ be a 3-uniform hypergraph with $n$ vertices and at least $Cn^{8/3}$ edges. Fix $t\geq 4$ and set $p=t^{-3}$ and $\eps=2t^{-6}$. We show that if $C$ is sufficiently large in terms of $t$ then $H$ contains a homeomorph of $K_t^{(3)}$.

\smallskip

\textit{Notation.} Let $G\subseteq\skel(H)$ be a graph. Given vertices $v_1,\ldots,v_k\in V(G)$, let $N_G(v_1,\ldots,v_k)$ denote their common neighborhood in $G$ and write $\codeg_G(v_1,\ldots,v_k)=|N_G(v_1,\ldots,v_k)|$ for the size of this set. Given $v,v'\in V(G)$, let $\xi_G(v,v')$ count the number of non-$(p,\eps)$-disk-coverable 4-cycles $vwv'w'\subset G$ that contain $v$ and $v'$ as opposite vertices. Define
\[
\psi_G(v,v')=\frac{\xi_G(v,v')}{\codeg_G(v,v')}
,\]
where the quotient is considered to be 0 if $\codeg_G(v,v')=0$.

\smallskip

The proof proceeds in three steps. In \cref{claim:findG}, we identify a subgraph $G$ of $\skel(H)$ in which $\psi_G$ is not too large on average. Then, in \cref{claim:findvi,claim:findGamma}, we use the dependent random choice technique to embed $\Gamma_t$ in $G$ such that every special 4-cycle of $\Gamma_t$ is $(p,\eps)$-disk-coverable in the embedding. Lastly, we use \cref{lem:embedding} to attach disks to the special 4-cycles of $\Gamma_t$.

We remark that $\psi_G(v,v')$ interpolates between two more natural quantities, $\xi_G(v,v')$ and $\xi_G(v,v')/\codeg(v,v')^2$; the latter measures the probability that a randomly chosen 4-cycle $vwv'w'$ is not disk-coverable. However, it proves most natural to bound $\psi_G(v,v')$ in the following claim. Additionally, controlling $\psi_G(v,v')$ later allows us to bound the probability that a randomly chosen 4-cycle $vwv'w'$ is disk-coverable, even if the choice of $w'$ is restricted to some small subset of $N_G(v_,v')$.

\begin{claim}\label{claim:findG}
There is a graph $G\subseteq\skel(H)$ with $e(G)> Cn^{5/3}$ and
\[
\sum_{\{v,v'\}\subset V(G)}\psi_G(v,v')<\frac{n^{1/3}}{2Cp^2\eps}e(G),
\]
where the sum is taken over unordered pairs of distinct vertices $\{v,v'\}\subset V(G)$.
\end{claim}

\begin{proof}\proofofclaim
Choose $u\in V(H)$ uniformly at random. We show that taking $G=H_u$ satisfies the given conditions with positive probability. First, we have
\[
\E[e(H_u)]=\sum_{e\in E(H)}\Pr[u\in e]=\frac 3ne(H)\geq 3Cn^{5/3}.
\]
Next, write $\Psi_u=\sum_{\{v,v'\}\subset V(H_u)}\psi_{H_u}(v,v')$. We have
\begin{align*}
\E\left[\Psi_u\right]
=\frac 1n\sum_{\substack{u\in V(H),\\\{v,v'\}\subset V(H_u)}}
	\frac{|\{\text{4-cycles $vwv'w'\subset H_u$ : $vwv'w'$ is not $(p,\eps)$-disk-coverable}\}|}{\codeg_{H_u}(v,v')},
\end{align*}
where the summand is defined to be 0 if $\codeg_{H_u}(v,v')=0$.

For any two vertices $v,v'\in V(H)$, write $\LL(v,v')=H_v\cap H_{v'}$. We have
\[N_{H_u}(v,v')=\{w\in V(H):uvw,uv'w\in E(H)\}=N_{\LL(v,v')}(u),\]
so $\codeg_{H_u}(v,v')=\deg_{\LL(v,v')}(u)$. Additionally, \cref{prop:diskc} implies that if $vwv'w'\subset H_u$ is not $(p,\eps)$-disk-coverable then the path $wuw'$ is not $(p,\eps)$-admissible in $\LL(v,v')$. Thus,
\[
\E[\Psi_u]\leq\frac 1n
	\sum_{\substack{\{v,v'\}\subset V(H),\\u\in V(\LL(v,v'))}}
	\frac{|\{\text{$wuw'\subset\LL(v,v')$ : $wuw'$ is not $(p,\eps)$-admissible}\}|}{\deg_{\LL(v,v')}(u)},
\]
where the summand is again 0 if $\deg_{\LL(v,v')}(u)=0$.
Applying \cref{lem:2adm} to each graph $\LL(v,v')$, we derive
\[
\E[\Psi_u]<\frac 1n\sum_{\{v,v'\}\subset V(H)}\frac{3|V(\LL(v,v'))|}{2p^2\eps}<\frac 1n\times \binom n2\times\frac{3n}{2p^2\eps}<\frac{3n^2}{4p^2\eps}.
\]

Combining our two expected value estimates,
\[
\E\left[n^{1/3}e(H_u) - Cn^2 - 2Cp^2\eps\Psi_u\right]>3Cn^2 - Cn^2 - \frac 32 Cn^2>0.
\]
Thus, there is $u\in V(H)$ such that $G=H_u$ satisfies $n^{1/3}e(G)-Cn^2-2Cp^2\eps\Psi_u>0$. This immediately implies that $e(G)>Cn^2/n^{1/3}=Cn^{5/3}$ and
\[\Psi_u=\sum_{\{v,v'\}\subset V(G)}\psi_G(v,v')<\frac{n^{1/3}}{2Cp^2\eps}e(G),
\]
as desired.
\end{proof}

Fix $G\subset\skel(H)$ as in \cref{claim:findG}; by adding isolated vertices to $G$, we may assume that $|V(G)|=|V(H)|=n$. The remainder of the proof focuses on identifying a suitable copy of $\Gamma_t$ in $G$. Henceforth, we drop the subscripts on $\xi$, $\psi$, $N$, and $\codeg$, which are always implicitly evaluated with respect to the graph $G$.

For vertices $v_1,v_2,v_3\in V(G)$, set
\[
\phi(v_1,v_2,v_3)=\frac{\psi(v_1,v_2)+\psi(v_1,v_3)+\psi(v_2,v_3)}{\codeg(v_1,v_2,v_3)}.
\]
The next claim identifies $t$ distinct vertices $v_1,\ldots,v_t\in V(G)$ such that each triple $\{v_i,v_j,v_k\}$ has many common neighbors and $\phi(v_i,v_j,v_k)$ is small on average. The vertices $v_1,\ldots,v_t$ constitute one side of an embedding of $\Gamma_t$ in $G$. The conditions on $v_1,\ldots,v_t$ guaranteed by \cref{claim:findvi} are then used in \cref{claim:findGamma} to embed the remaining $\binom t2+\binom t3$ vertices of $\Gamma_t$ such that every special 4-cycle of $\Gamma_t$ is $(p,\eps)$-disk-coverable.

\begin{claim}\label{claim:findvi}
Fix $r\geq 0$ satisfying $r+3/p^2\eps\leq C^3/t^3$. Then, there are $t$ distinct vertices $v_1,\ldots,v_t\in V(G)$ such that $\codeg(v_i,v_j,v_k)> r$ for any $i,j,k\in[t]$ and
\[
\sum_{1\leq i<j<k\leq t}\phi(v_i,v_j,v_k)<\frac 16.
\]	
\end{claim}

\begin{proof}\proofofclaim
We proceed in two steps. First, we identify $W\subset V(G)$ such that most triples in $W$ have large common neighborhood and small $\phi$. Then, we choose $v_1,\ldots,v_t\in W$.

The first step is a standard application of the dependent random choice method (see \cite{FoSu11} for further discussion of this technique). Choose $w\in V(G)$ uniformly at random and let $W=N(w)$. Observe that
\[\E[|W|]=2e(G)/n>2Cn^{2/3}.\]
Let $T\subset\binom W3$ be the set of triples $\{w_1,w_2,w_3\}$ in $W$ with $\codeg(w_1,w_2,w_3)\leq r$. We have
\begin{align*}
\E[|T|]&=\sum_{\substack{\{w_1,w_2,w_3\}\subset V(G),\\\codeg(w_1,w_2,w_3)\leq r}}\Pr[\{w_1,w_2,w_3\}\subset W]
=\sum_{\substack{\{w_1,w_2,w_3\}\subset V(G),\\\codeg(w_1,w_2,w_3)\leq r}}\frac{\codeg(w_1,w_2,w_3)}n
\\&\leq\binom n3 \frac rn<\frac{rn^2}6
<\frac{rn^{1/3}}{6C}e(G).
\end{align*}
Second, set 
\[
\Phi=\sum_{\{w_1,w_2,w_3\}\subset W}\phi(w_1,w_2,w_3).
\]
We have
\begin{align*}
\E[\Phi]&=\sum_{\{w_1,w_2,w_3\}\subset V(G)}\phi(w_1,w_2,w_3)\times\Pr[\{w_1,w_2,w_3\}\subset W]
\\&=\sum_{\{w_1,w_2,w_3\}\subset V(G)}
	\frac {\psi(w_1,w_2)+\psi(w_1,w_3)+\psi(w_2,w_3)}n
\\&=\frac{n-2}n\sum_{\{w_1,w_2\}\subset V(G)}\psi(w_1,w_2)
<\frac{n^{1/3}}{2Cp^2\eps}e(G),
\end{align*}
where the final inequality comes directly from \cref{claim:findG}.

Combining the prior bounds,
\begin{align*}
\E&\left[|W|-\frac{e(G)}n-\left(\frac{6C}{r+3/p^2\eps}\right)n^{-4/3}(|T|+\Phi)\right]
\\&\quad>\frac{2e(G)}n-\frac{e(G)}n-\left(\frac{6C}{r+3/p^2\eps}\right)\left(\frac{r}{6C}+\frac 1{2Cp^2\eps}\right)\frac{e(G)}n = 0.
\end{align*}
Thus, we may fix $W\subset V(G)$ with
\[
|W|-\frac{e(G)}n-\left(\frac{6C}{r+3p^2\eps}\right)n^{-4/3}(|T|+\Phi)>0.
\]
In particular, $|W|> e(G)/n> Cn^{2/3}$,
and
\[
|T|+\Phi<\left(\frac{r+3/p^2\eps}{6C}\right)n^{4/3}|W|<\left(\frac{r+3/p^2\eps}{6C^3}\right)|W|^3.
\]

Having fixed $W$, choose $v_1,\ldots,v_t\in W$ uniformly at random without replacement. This is permissible because $|W|>Cn^{2/3}\geq C> t$. Consider the random variable
\[
X=\sum_{1\leq i<j<k\leq t}\mathbf 1{\{\codeg(v_i,v_j,v_k)\leq r\}}+\phi(v_i,v_j,v_k).
\]
We have
\[
\E[X]=\frac{t(t-1)(t-2)}{|W|(|W|-1)(|W|-2)}\times (|T|+\Phi)
<\frac{t^3}{|W|^3}\times\left(\frac{r+3/p^2\eps}{6C^3}\right)|W|^3=\frac{t^3(r+3/p^2\eps)}{6C^3}\leq\frac 16.
\]
Thus, there is a choice of $v_1,\ldots,v_t\in W$ such that $X<1/6$. In particular, no trio $\{v_i,v_j,v_k\}$ satisfies $\codeg(v_i,v_j,v_k)\leq r$, and
\[
\sum_{1\leq i<j<k\leq t}\phi(v_i,v_j,v_k)\leq X<\frac 16.\qedhere
\]
\end{proof}

We now embed the remainder of $\Gamma_t$ in $G$.

\begin{claim}\label{claim:findGamma}
Suppose $C\geq\left(t^9+3t^3/p^2\eps\right)^{1/3}$. 
Then there is an embedding of $\Gamma_t$ in $G$ such that every special 4-cycle of $\Gamma_t$ is $(p,\eps)$-disk-coverable.
\end{claim}

\begin{proof}\proofofclaim
Observe that $C$ satisfies the hypothesis of \cref{claim:findvi} with $r=t^6$. Fix vertices $v_1,\ldots,v_t\in V(G)$ satisfying \cref{claim:findvi} for this value of $r$.
For each $ij\in\binom{[t]}2$, choose $v_{ij}\in N(v_i,v_j)$ uniformly at random. Similarly, for each $ijk\in\binom{[t]}3$, choose $v_{ijk}\in N(v_i,v_j,v_k)$ uniformly at random.

Consider the event that this is an embedding, i.e., that the vertices $v_S$ are pairwise distinct for all $S\in\binom{[t]}1\cup\binom{[t]}2\cup\binom{[t]}3$. If $S,T\in\binom{[t]}1\cup\binom{[t]}2\cup\binom{[t]}3$ do not both have size 1 then $\Pr[v_S=v_T]< t^{-6}$, because each vertex $v_{ij}$ or $v_{ijk}$ is chosen uniformly among more than $r=t^6$ possibilities. Recalling that the vertices $v_1,\ldots,v_t$ are pairwise distinct, it follows that the expected number of equalities among the $\binom t1+\binom t2+\binom t3$ vertices $v_S$ is at most
\[
\left[\frac 12\left(\binom t2+\binom t3\right)^2+\binom t1\left(\binom t2+\binom t3\right)\right]\times\frac 1{t^6}
<\left[\frac 12\left(\frac{t^3}6\right)^2+t\times\frac{t^3}6\right]\times\frac 1{t^6}
=\frac 1{72}+\frac 1{6t^2}
<\frac 1{36}.
\]
Thus, with probability at least $35/36$, the vertices $v_S$ are pairwise distinct for all $S\in\binom{[t]}1\cup\binom{[t]}2\cup\binom{[t]}3$, yielding a nondegenerate embedding of $\Gamma_t$ in $G$.

Next, we consider the event that the embedding of some special 4-cycle $v_iv_{ij}v_jv_{ijk}$ is nondegenerate but not $(p,\eps)$-disk-coverable. Given $ij\in\binom{[t]}2$ and $k\in [t]\setminus\{i,j\}$, let $A_{ij;k}$ be the event that $v_{ij}\neq v_{ijk}$ and the 4-cycle $v_iv_{ij}v_jv_{ijk}$ is not $(p,\eps)$-disk-coverable. Note that
\[
\Pr[A_{ij;k}]\leq\frac{2\xi(v_i,v_j)}{\codeg(v_i,v_j)\codeg(v_i,v_j,v_k)}=\frac{2\psi(v_i,v_j)}{\codeg(v_i,v_j,v_k)}.
\]
Using a union bound and then applying \cref{claim:findvi} for the final inequality, we have
\begin{align*}
\Pr[A_{ij;k}\text{ holds for some }i,j,k]
&\leq\sum_{1\leq i<j<k\leq t}\Pr[A_{ij;k}]+\Pr[A_{ik;j}]+\Pr[A_{jk;i}]
\\&\leq\sum_{1\leq i<j<k\leq t}\frac{2(\psi(v_i,v_j)+\psi(v_i,v_k)+\psi(v_j,v_k))}{\codeg(v_i,v_j,v_k)}
\\&=\sum_{1\leq i<j<k\leq t}2\phi(v_i,v_j,v_k)<\frac 13.
\end{align*}
It follows that with probability at least $1-\frac 1{36}-\frac 13>0$, the vertices $v_S$ are pairwise distinct for all $S\in\binom{[t]}1+\binom{[t]}2+\binom{[t]}3$, and furthermore, each special 4-cycle $v_iv_{ij}v_jv_{ijk}$ is $(p,\eps)$-disk-coverable.	
\end{proof}

Lastly, we attach disks to the special 4-cycles of $\Gamma_t$, leveraging their $(p,\eps)$-disk-coverability. The number of special 4-cycles $v_iv_{ij}v_jv_{ijk}$ in $\Gamma_t$ is $k=3\binom t3<t^3/2$ and its vertex set has size $m=\binom t1+\binom t2+\binom t3<t^3/2$. The parameters $p=t^{-3}$ and $\eps=2t^{-6}$ satisfy $p<\frac 1{k+m}$ and $\eps<\frac 1{k(k+m)}$, so \cref{lem:embedding} applies. This yields $3\binom t3$ disks, each bounded by one of the special 4-cycles of $\Gamma_t$, such that any two disks (when viewed as simplicial complexes) intersect only on their boundaries. It follows that the union of these $3\binom t3$ disks is a homeomorph of $K_t^{(3)}$.

\section{The Torus and the Real Projective Plane}\label{sec:TandRP}

In this section, we leverage \cref{lem:2adm} to give a shorter proof of \cref{thm:surfaces}.
We begin with a lemma that is a direct corollary of results from \cite{KPTZ21} and \cite{Sa22}.

\begin{lemma}%[cf.\ {\cite[Lemma 3.3]{KPTZ21}} + {\cite[Lemma 3.3]{Sa22}}]
\label{lem:P1adm}
For all $p,\eps\in(0,1)$ there is a constant $C=C(p,\eps)$ such that the following holds. Let $H$ be a 3-uniform hypergraph with at least $Cn^{5/2}$ edges. Then, there is a subset $F\subseteq E(H)$ of size $|F|\geq e(H)/2$ such that, for any neighboring edges $xyz,x'yz\in F$, the 4-cycle $xyx'z$ bounding their union is $(p,\eps)$-disk-coverable in $H$.
\end{lemma}

\begin{proof}
Choose an integer $r$ large enough that $(1-p)^r<\eps/2$ and set $\eps'=\eps/4r$. Lemma 3.3 in \cite{KPTZ21} states that if $e(H)\geq\frac{24r}{p}\sqrt\frac 2{\eps'}$ then there is a subset $F\subseteq E(H)$ containing at least $e(H)/2$ edges such that any neighboring edges $xyz,x'yz\in F$ are \emph{$(p/2,\eps',2,r)$-semi-admissible}, a property defined in \cite{KPTZ21}. The definition of semi-admissibility is somewhat involved, so we do not give it here. Instead, we cite Lemma 3.3 from \cite{Sa22}, which states that any $(p/2,\eps',2,r)$-semi-admissible pair $(xyz,x'yz)$ gives rise to a $(p,2r\eps'+(1-p)^r)$-disk-coverable 4-cycle $xyx'z$. 
 Our proof concludes with the observation that $2r\eps'+(1-p)^r<\eps$, by the choice of $r$ and $\eps'$.
\end{proof}

We now prove \cref{thm:surfaces}, which is restated in the following equivalent form.

\begin{theorem}
	There is a constant $C$ such that the following holds. Any 3-uniform hypergraph $H$ with at least $Cn^{5/2}$ edges contains homeomorphs of both the torus $\TT$ and the real projective plane $\RP$.
\end{theorem}

\begin{proof}
Set $p=1/18$ and $\eps=1/163$. Suppose $C$ is a sufficiently large constant --- more precisely, we require that the conclusion of \cref{lem:P1adm} holds for $p,\eps$ and that $6C^2\geq 30p^2\eps + 15$.

Let $H$ be an $n$-vertex 3-uniform hypergraph with at least $Cn^{5/2}$ edges and let $V$ be its vertex set. Applying \cref{lem:P1adm}, we obtain a set $F\subseteq E(H)$ of at least $Cn^{5/2}/2$ edges such that $xyx'z$ is $(p,\eps)$-disk-admissible in $H$ for each pair of neighboring edges $xyz,x'yz\in F$. Let $H'$ be the 3-uniform hypergraph on $V$ with edge set $E(H')=F$.

\begin{claim}
	There are distinct vertices $u_1,u_2\in V$ such that the graph $H'_{u_1}\cap H'_{u_2}$ has at least $3C^2n$ edges.
\end{claim}

\begin{proof}[Proof of Claim]\proofofclaim
Choose $u_1,u_2\in V$ uniformly at random without replacement. Given an unordered pair of vertices $vw\in \binom V2$, set $T_{vw}=|\{u\in V:uvw\in F\}|$. Then,
\[
\E_{u_1,u_2}\left[e(H'_{u_1}\cap H'_{u_2})\right]
=\frac 1{n(n-1)}\sum_{u_1\neq u_2}e(H'_{u_1}\cap H'_{u_2})=\frac 1{n(n-1)}\sum_{vw\in\binom V2}T_{vw}(T_{vw}-1).
\]
The average of the $T_{vw}$ is
\[
T=\frac 2{n(n-1)}\sum_{vw\in\binom V2}T_{vw}=\frac 2{n(n-1)}\times 3|F|\geq 3Cn^{1/2}.
\]
By convexity, it follows that
\begin{align*}
\E_{u_1,u_2}\left[e(H'_{u_1}\cap H'_{u_2})\right]
&=\frac 1{n(n-1)}\times \frac{n(n-1)}{2}T(T-1)\geq\frac 12(3Cn^{1/2})(3Cn^{1/2}-1)
\\&\geq\frac 12(3Cn^{1/2})(2Cn^{1/2})=3C^2n.
\end{align*}
Thus there is a choice of $u_1\neq u_2$ such that $H'_{u_1}\cap H'_{u_2}$ has at least $3C^2n$ edges.
\end{proof}

Fix $u_1,u_2\in V$ such that the graph $G=H'_{u_1}\cap H'_{u_2}$ has at least $3C^2n$ edges.

\begin{claim}\label{claim:TRP-2}
	There is a vertex $v\in V(G)$ with six distinct neighbors $w_1,\ldots w_6\in N_G(v)$ such that the three paths $w_1vw_2$, $w_3vw_4$, and $w_5vw_6$ are each $(p,\eps)$-admissible in $G$.
\end{claim}

\begin{proof}[Proof of Claim]\proofofclaim
For each $v\in V(G)$, let $b_v$ be the number of ``bad'' length-2 paths $wvw'$ in $G$ that are not $(p,\eps)$-admissible. 
We claim that there is some vertex $v\in V(G)$ with $\deg(v)\geq 15$ and $b_v\leq(\deg v)^2/20$. Observe that
\[
\sum_{\substack{v\in V(G),\\\deg v<15}}\deg v + \sum_{\substack{v\in V(G),\\b_v>(\deg v)^2/20}}\deg v
< 15n + 20\sum_{v\in V(G)}\frac{b_v}{\deg v}\leq 15n+\frac{30n}{p^2\eps},
\]
by \cref{lem:2adm}. Recalling that $6C^2\geq 30/p^2\eps + 15$, we have
\[
\sum_{\substack{v\in V(G),\\\deg v<15}}\deg v + \sum_{\substack{v\in V(G),\\b_v>(\deg v)^2/20}}\deg v
<6C^2n\leq 2e(G)=\sum_{v\in V(G)}\deg v.
\]
Thus there is a vertex $v\in V(G)$ contributing to neither summand on the left, i.e., $v$ satisfies $\deg(v)\geq 15$ and $b_v\leq(\deg v)^2/20$.

Choose $w_1,\ldots,w_6\in N_G(v)$ uniformly at random with replacement. Observe that
\[
\Pr[w_1,\ldots,w_6\text{ distinct}]=\frac{15\times 14\times\cdots\times 10}{15^6}\geq 0.31
\]
and for each $i=1,3,5$,
\[
\Pr[\text{$w_i\neq w_{i+1}$ and $w_ivw_{i+1}$ not $(p,\eps)$-admissible}]=\frac{2b_v}{(\deg v)^2}\leq\frac 1{10}.
\]
Thus, a union bound gives
\[
\Pr[\text{$w_1,\ldots,w_6$ distinct and $w_ivw_{i+1}$ is $(p,\eps)$-admissible for $i=1,3,5$}]\geq 0.31-\frac 3{10}>0.
\]
It follows that there are distinct vertices $w_1,\ldots,w_6$ as desired.
\end{proof}

Fix vertices $v,w_1,\ldots,w_6$ as in \cref{claim:TRP-2}. By \cref{prop:diskc}, the 4-cycles $u_1w_iu_2w_{i+1}$ are $(p,\eps)$-disk-coverable in $H'$ for $i=1,3,5$. Moreover, for each $i=1,2$ and each $j\neq k$, the 4-cycle $u_iw_jvw_k$ is $(p,\eps)$-disk-coverable because $u_ivw_j$ and $u_ivw_k$ are neighboring edges in $F$.
%Moreover, for each $i\neq j$, the 4-cycles $uw_ivw_j$ and $u'w_ivw_j$ are $(p,\eps)$-disk-coverable because $(uvw_i,uvw_j)$ and $(u'vw_i,u'vw_j)$ are pairs of neighboring edges in $F$.

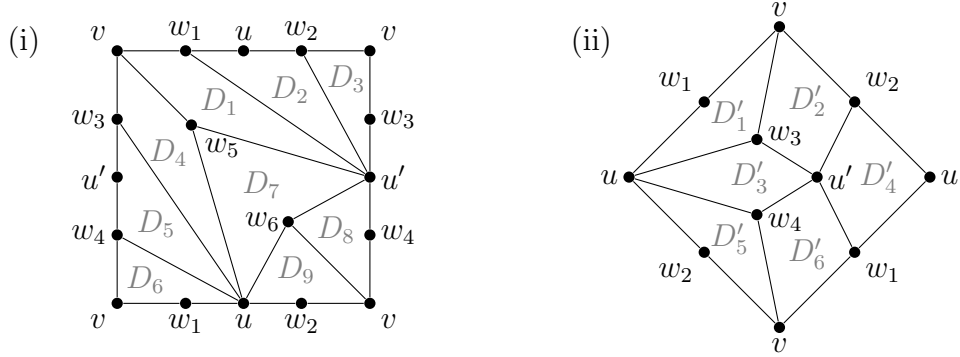
\begin{figure}[t]
\centering
\begin{tikzpicture}[baseline={(0,0)}]
\node at (0,1.8) {(i)};	
\end{tikzpicture}
\begin{tikzpicture}[scale=1.4, baseline={(0,0)}]
	%\node at (-1.9,2/1.4) {(ii)};
	\foreach \ang/\pos [count=\i] in {135/above left,45/above right,-45/below right,-135/below left} {
		\coordinate[vtx] (v0\i) at (\ang:1.414*1.2);
		\node[\pos] at (v0\i) {$v$};
	}
	\foreach \lab/\nam [count=\i, evaluate=\x using {(\i - 2) * 0.55}]
	in {v1/w_1,u1/u_1,v2/w_2} {
		\coordinate[vtx] (\lab1) at (\x,1.2);
		\node[above] at (\lab1) {$\nam$};
		\coordinate[vtx] (\lab2) at (\x,-1.2);
		\node[below] at (\lab2) {$\nam$};
	}
	\foreach \lab/\nam [count=\i, evaluate=\y using {(\i - 2) * 0.55}]
	in {v4/w_4,u2/u_2,v3/w_3} {
		\coordinate[vtx] (\lab1) at (-1.2,\y);
		\node[left] at (\lab1) {$\nam$};
		\coordinate[vtx] (\lab2) at (1.2,\y);
		\node[right] at (\lab2) {$\nam$};
	}
	\draw (v01) -- (v11) -- (u11) -- (v21)
	   -- (v02) -- (v32) -- (u22) -- (v42)
	   -- (v03) -- (v22) -- (u12) -- (v12)
	   -- (v04) -- (v41) -- (u21) -- (v31) -- (v01);
	\path (135:0.7) coordinate[vtx] (v5) node[below right=1] {$w_5$}
		  (-45:0.6) coordinate[vtx] (v6) node[left=-1] {$w_6$};
	\draw (v11) -- (u22) -- (v21)
	      (v31) -- (u12) -- (v41)
	      (v01) -- (v5)
	      (v6) -- (v03)
	      (v5) -- (u12) -- (v6) -- (u22) -- (v5);
	\path[gray] (v01) -- (u22) node[pos=0.4] {$D_1$}
	      (v01) -- (u12) node[pos=0.4] {$D_4$}
	      (u11) --++(0.1,0) -- (u22) node[pos=0.32] {$D_2$}
	      (u21) --++(0,-0.1) -- (u12) node[pos=0.32] {$D_5$}
	      (v02) --++(-130:0.35) node {$D_3$}
	      (v04) --++(40:0.35) node {$D_6$}
	      (-45:0.1) --++(0.1,0) node{$D_7$}
	      (v6) -- (v42) node[pos=0.6] {$D_8$}
	      (v6) -- (v22) node[pos=0.6] {$D_9$};
\end{tikzpicture}
\qquad\qquad
\begin{tikzpicture}[baseline={(0,0)}]
	\node at (-2.5,1.8){(ii)};
	\path (-2,0) coordinate[vtx] (u1) node[left] {$u_1$}
	   -- (-1,1) coordinate[vtx] (v1) node[above left] {$w_1$}
	   -- (0,2) coordinate[vtx] (v0) node[above] {$v$}
	   -- (1,1) coordinate[vtx] (v2) node[above right] {$w_2$}
	   -- (2,0) coordinate[vtx] (u2) node[right] {$u_1$}
	   -- (1,-1) coordinate[vtx] (w1) node[below right] {$w_1$}
	   -- (0,-2) coordinate[vtx] (w0) node[below] {$v$}
	   -- (-1,-1) coordinate[vtx] (w2) node[below left] {$w_2$}
	   -- (-0.3,0.5) coordinate[vtx] (v3) --+ (0,0.05) node[right] {$w_3$}
	   -- (-0.3,-0.5) coordinate[vtx] (v4) --+ (0,-0.1) node[right] {$w_4$}
	   -- (0.5,0) coordinate[vtx] (u3) --++ (-28:11pt) node {$u_2$};
	\draw (u1) -- (v1) -- (v0) -- (v2) -- (u2) -- (w1) -- (w0) -- (w2) -- (u1);
	\draw (u1) -- (v3) -- (u3) -- (v4) -- (u1)
	      (v0) -- (v3)
	      (w0) -- (v4)
	      (v2) -- (u3) -- (w1);
	\path[gray] (u1) -- (u3) node[pos=0.65] {$D'_3$}
	      (u3) -- (u2) node[above=-7pt,pos=0.55] {$D'_4$}
	      (v3) node[above left=-1] {$D'_1$}
	      (v4) node[below left=-1] {$D'_5$}
	      (70:1.08) node {$D'_2$}
	      (-70:1.08) node {$D'_6$};
\end{tikzpicture}
\caption{Decompositions of (i) the torus and (ii) the real projective plane as CW complexes.}
\label{fig:quadrangulations}
\end{figure}

We construct homeomorphs of $\TT$ and $\RP$ by using \cref{lem:embedding} to attach disks to those $(p,\eps)$-disk-coverable 4-cycles pictured in \cref{fig:quadrangulations}. For the torus, it suffices to attach disks $D_1,\ldots,D_9$ to the nine $(p,\eps)$-disk-coverable 4-cycles
\[
\begin{array}{ccccc}
C_1=u_2w_1vw_5,
&C_2=u_1w_1u_2w_2,
&C_3=u_2w_2vw_3,
&C_4=u_1w_3vw_5, 
&C_5=u_1w_3u_2w_4,
\\&C_6=u_1w_1vw_4, 
&C_7=u_1w_5u_2w_6,
&C_8=u_2w_4vw_6, 
&C_9=u_1w_2vw_6
\end{array}
\]
pictured in \cref{fig:quadrangulations}(i). We may apply \cref{lem:embedding} with $k=m=9$ to locate suitable disks, because $p=\frac 1{9+9}$ and $\eps<\frac 1{9(9+9)}$.
For the real projective plane, it suffices to attach disks $D'_1,\ldots,D'_6$ to the six $(p,\eps)$-disk-coverable 4-cycles
\[
\begin{array}{cccccc}
C'_1=u_1w_1vw_3,
&C'_2=u_2w_2vw_3,
&C'_3=u_1w_3u_2w_4,
\\C'_4=u_1w_1u_2w_2,
&C'_5=u_1w_2vw_4,
&C'_6=u_2w_1vw_4
\end{array}
\]
pictured in \cref{fig:quadrangulations}(ii). We apply \cref{lem:embedding} again, this time with parameters $k=6$ and $m=7$.
\end{proof}

\section{Concluding Remarks}\label{sec:end}

In this paper, we showed that $\exh(n,X)=O(n^{8/3})$ for any fixed 2-dimensional simplicial complex $X$. Like the authors in \cite{KLNS21}, we suspect that the correct exponent is $5/2$, i.e., that \cref{conj:upperbound} holds.

This section highlights connections between our approach and the study of bipartite Tur\'an numbers --- including some limitations imposed thereby. We then conclude with a brief discussion of lower bounds in the 2-dimensional case and open questions in higher dimensions.

\subsection{Connection to bipartite Tur\'an numbers}
It is not clear whether our method can be pushed to achieve a better exponent for $\exh(n,K_t^{(3)})$, because the asymptotics of the standard Tur\'an number $\ex(n,\Gamma_t)$ remain unknown. In particular, one cannot necessarily embed $\Gamma_t$ in a link graph comprising fewer than $\ex(n,\Gamma_t)$ edges. Thus, our approach is inherently limited to 3-uniform hypergraphs $H$ with $n$ vertices and at least $\Omega(n\times \ex(n,\Gamma_t))$ edges.

Until recently, the best-known upper bound for $\ex(n,\Gamma_t)$ was $O(n^{5/3})$, proven by F\"uredi \cite{Fu91} in 1991 and reproved in 2003 by Alon, Krivelevich, and Sudakov \cite{AlKrSu03} using dependent random choice. In fact, the latter argument may be recovered from \cref{claim:findvi,claim:findGamma} if one ignores the disk-coverability constraint. The upper bound was improved to $o(n^{5/3})$ by Sudakov and Tomon \cite{SuTo20} in 2020, but it remains open to improve this by a polynomial factor. Thus, any improvement to the topological Tur\'an exponent of $K_t^{(3)}$ by our approach would first require improving the Tur\'an exponent of $\Gamma_t$.

In this spirit, we raise the following question regarding the asymptotics of $\ex(n,\Gamma_t)$.

\begin{question}\label{q:exGamma}
Is $\ex(n,\Gamma_t)=O(n^{5/3-\eps})$ for some $\eps>0$, which depends on $t$?
\end{question}

\noindent
By a simple random construction, the Tur\'an exponent of $\Gamma_t$ is at least $2-v(\Gamma_t)/e(\Gamma_t)=5/3-O(1/t)$. Hence, an affirmative answer to \cref{q:exGamma} could improve the upper bound on $\exh(n,K_t^{(3)})$ by a factor of $n^{O(1/t)}$ at best. However, an affirmative answer to \cref{q:exGamma} would provide evidence towards our conjecture that there is a better universal exponent, particularly if we could bootstrap the proof to find an embedding of $\Gamma_t$ with disk-coverable special 4-cycles in the 1-skeleton of any 3-uniform hypergraph with $\Omega(n^{8/3-\eps})$ edges.

At this point, it is natural to search for an alternative CW complex representation of $K_3^{(t)}$ whose 1-skeleton is known to have Tur\'an exponent less than 5/3. This is not difficult if we permit disks to be attached to 6-cycles instead of 4-cycles. Consider for example the subgraph $\Gamma'_t\subseteq\Gamma_t$ induced by the $\binom t1+\binom t2$ vertices of the form $v_i$ or $v_{ij}$. The upper bound from \cite{AlKrSu03,Fu91} yields $\ex(n,\Gamma'_t)=O(n^{3/2})$, and attaching a disk to each 6-cycle of the form $v_iv_{ij}v_jv_{jk}v_kv_{ki}$ results in a homeomorph of $K_t^{(3)}$. However, there is no obvious way to alter the proof of \cref{claim:findG} --- which bounds the number of non-disk-coverable 4-cycles in a random link $H_u$ --- to derive an analogous result for 6-cycles.

% despite recent progress OF AUTHORS [], this remains

%However, we were unable to count non-disk-coverable 6-cycles in a manner analogous to \cref{claim:findG}; we leave this as an open problem.
%
%\begin{problem}
%In \cref{claim:findG}, we bound the expected number of non-disk-coverable 4-cycles in a random link $H_u$, summed with appropriate weighting factors --- this is $\E[\Psi]$. Is there an analogous result for non-disk-coverable 6-cycles?
%\end{problem}

If we instead choose a CW representation of $K_t^{(3)}$ in which disks are only attached to 4-cycles, then it is not difficult to check that its 1-skeleton $\Gamma$ cannot be too sparse. In particular, if $t\geq 5$, then $\Gamma$ cannot be 2-degenerate; equivalently, it must contain a subgraph of minimum degree 3. Thus, we cannot identify a suitable CW representation of $K_3^{(t)}$ without first substantially improving our understanding of the Tur\'an exponents of such graphs. Despite some recent progress \cite{Ja23,SuTo20}, this remains a major open problem.

\subsection{Lower Bounds}

Keevash, Polyanskii, Tomon, and Zakharov \cite{KPTZ21} prove a lower bound matching \cref{conj:upperbound} when $X$ is a surface, showing that $\exh(n,X)\geq cn^{5/2}$ for some positive constant $c=c(X)$. We conjecture that this lower bound extends to all sufficiently complicated $X$.

\begin{conjecture}\label{conj:lowerbound}
Let $X$ be a 2-dimensional simplicial complex in which every 1-simplex is contained in at least two 2-simplices. Then there is a positive constant $c=c(X)$ such that $\exh(n,X)\geq cn^{5/2}$.
\end{conjecture} % there are spaces like RP2 that satisfy ^ but do not have nontrivial second homology group

The lower bound in \cite{KPTZ21} follows from a result of Gao \cite{Ga91} showing that the number of $n$-vertex triangulations of a fixed surface $X$ is exponential in $n$. A natural approach towards proving \cref{conj:lowerbound} would be to show a similar result for any 2-dimensional simplicial complex $X$. Such a result would immediately yield \cref{conj:lowerbound}, by the argument from \cite{KPTZ21}.

We observe that \cref{conj:lowerbound} is false if even a few 1-simplices of $X$ are not contained in multiple 2-simplices. For example, deleting a 2-simplex from any triangulation of the sphere $\mathbb S^2$ yields a triangulation of the disk $\disk$ in which all but three 1-simplices are contained in multiple 2-simplices. However, it is clear that $\exh(n,\disk)=0$.

\subsection{Higher Dimensions}

In dimensions $d\geq 3$, very little is known about topological Tur\'an exponents. For spheres, the best-known upper bound on $\exh(n,\mathbb S^d)$ is $O(n^{d+1-2^{1-d}})$ due to Newman and Pavelka \cite{NePa24}, which follows from a short Cauchy--Schwarz argument. Newman and Pavelka conjecture that this exponent can be improved.
\begin{conjecture}[{\cite[Conjecture 13]{NePa24}}]\label{conj:sph}
	It holds that $\exh(n,\mathbb S^d)=O(n^{d+1-(d+1)/(2^{d+1}-2)})$.
\end{conjecture}
\noindent
Indeed, any improvement to the Cauchy--Schwarz upper bound would be quite interesting.

We remark that proving any lower bound on $\exh(n,\mathbb S^d)$ is quite challenging, as triangulations of $\mathbb S^d$ do not necessarily admit any nice combinatorial structure in general. Newman and Pavelka showed that \cref{conj:sph} is tight if we restrict our focus to various structured families of triangulations --- such as shellable spheres --- and is potentially tight in general, assuming an affirmative answer to an open problem of Gromov.

Lastly, it would be interesting to improve our understanding of the universal exponent in higher dimensions.
\begin{question}
What is the largest constant $\eps_d$ such that, for any $d$-dimensional simplicial complex $X$, we have $\exh(n,X)=n^{d+1-\eps_d+o(1)}$?
\end{question}

As mentioned in the introduction, Long, Narayanan, and Yap \cite{LNY22} showed that $\eps_d\geq d^{-2d^2}$; however, there is no reason to believe that this is even asymptotically optimal. In recent upper bounds on 2-dimensional topological Tur\'an numbers \cite{KLNS21,KPTZ21,Sa22}, the strategy has relied on being able to find several copies of the sphere $\mathbb S^2$ which glue together in an appropriate configuration. Following this heuristic, it is not unreasonable to suspect that the universal exponent in higher dimensions behaves similarly to the exponent for the $d$-dimensional sphere.

\begin{conjecture}\label{conj:eps-exponential}
The universal exponent $\eps_d$ satisfies $\eps_d=2^{-d-o(d)}$.
\end{conjecture}

An interesting first step towards \cref{conj:eps-exponential} would be to show that $\eps_d$ is exponential in $d$, i.e., $\eps_d=2^{\Theta(-d)}$.

\end{document}